\documentclass[12pt,a4paper]{article}
\usepackage[cp1251]{inputenc}
\usepackage[english,ukrainian]{babel}
\usepackage{amsfonts,amssymb,amsmath}

\def\thebibliography#1{\centerline{\bf Література}
 \list
 {[\arabic{enumi}]}{\settowidth\labelwidth{[#1]}\leftmargin\labelwidth
 \advance\leftmargin\labelsep
 \usecounter{enumi}}
 \def\newblock{\hskip .11em plus .33em minus .07em}
 \sloppy\clubpenalty4000\widowpenalty4000
 \sfcode`\.=1000\relax}

\begin{document}

УДК 517.5

\setlength{\baselineskip}{16pt}

\setcounter{footnote}{0}

\begin{center}

{\renewcommand{\thefootnote}{ }

{\Large Один кубічний $3$-монотонний сплайн}

{\Large One cubic 3-monotone spline}

\vspace{0.2cm}

{\large Г.А.Дзюбенко}\\
{\large G.A.Dzyubenko}\\
{\small
{\it dzyuben@gmail.com}\\
Міжнародний математичний центр
 ім. Ю.О.Митропольського НАН України
 }
  }

\vspace{.8 cm}

\begin{minipage}{15 cm}

{\bf Abstract.}
For any 3-monotone on $[a,b]$ function $f$ (its third divided differences are nonnegative
for all choices of four distinct points, or equivalently, $f$ has a convex derivative on $(a,b)$) we construct a cubic 3-monotone (like $f$) spline $s$ with $n\in \Bbb N$ "almost" \ equidistant knots $a_j$ such that
$$
\left\Vert f-s \right\Vert_{[a_j,a_{j-1}]} \le c\, \omega_4
\left(f,(b-a)/n,[a_{j+4},a_{j-5}]\cap [a,b]\right), \quad
j=1,...,n,
$$
 where $c$ is an absolute constant, $\omega_4 \left(f,t,[\cdot,\cdot]\right)$
is the $4$-th modulus of smoothness of $f$,
and $||\cdot ||_{[\cdot,\cdot]}$ is the max-norm.

\vspace{.8 cm}

{\bf Анотація.}
Для будь-якої 3-монотонної на $[a,b]$ функції $f$ (її третя розділена різниця невід'ємна для всіх наборів з чотирьох різних точок, або еквівалентно, $f$ має опуклу на $(a,b)$ похідну), побудувано кубічний 3-монотонний (як і $f$) сплайн $s$ з $n\in \Bbb N$
"майже"\  рівновіддаленими вузлами $a_j$ такий, що
$$
\left\Vert f-s \right\Vert_{[a_j,a_{j-1}]} \le c\, \omega_4
\left(f,(b-a)/n,[a_{j+4},a_{j-5}]\cap [a,b]\right), \quad
j=1,...,n,
$$
 де $c $ --  абсолютна стала, $\omega_4 \left(f,t,[\cdot,\cdot]\right)$ --
 $4$-й модуль гладкості $f$
 і $||\cdot ||_{[\cdot,\cdot]} $ -- рівномірна норма.

\vspace{.8 cm}

{\bf Аннотация.}
Для любой 3-монотонной на $[a,b]$ функции $f$ (ее третья разделенная разность неотрицательна для всех наборов из четырех разных точек, или эквивалентно, $f$ имеет выпуклую на $(a,b)$ производную), построен кубический 3-монотонный (как $f$) сплайн $s$ с $n\in \Bbb N$ "почти"\ равноудаленными узлами $a_j$ такой, что
$$
\left\Vert f-s \right\Vert_{[a_j,a_{j-1}]} \le c\, \omega_4
\left(f,(b-a)/n,[a_{j+4},a_{j-5}]\cap [a,b]\right), \quad
j=1,...,n,
$$
 где $c$ -- абсолютная постоянная, $\omega_4 \left(f,t,[\cdot,\cdot]\right)$ --
 $4$-й модуль гладкости $f$
 и $||\cdot ||_{[\cdot,\cdot]} $ -- равномерная норма.

\end{minipage}
\end{center}

\bigskip
\medskip

\vspace{.3 cm}

\centerline{\bf \S 1. Вступ}
\vglue 0.3cm

Нахай $C:= C[a,b] $ -- простір неперервних на $[a,b]$
функцій $f:\, [a,b]\rightarrow \Bbb R$ з рівномірною нормою $ \|
f\|:=\| f\|_{[a,b]}:=\max\limits_{x\in [a,b]}|f(x)|,\
C^q:=\{f:\,f^{(q)}\in C\},\ q\in\Bbb N,$ і нехай  $\Delta^{3} :=\Delta^{3} [a,b]$ -- множина функцій
$f\in C,$ що мають невід'ємну третю розділену різницю в усіх наборах з чотирьох різних точок.
Зауважимо, що якщо $f\in\Delta^{3},$ то $f\in C^1$ і $f'$ опукла на $(a,b)$. Функції з
 $\Delta^{3}$ називаються {\it 3-монотонними.} Також, якщо $f\in C$ є тричі неперервно диференційованою на $(a,b),$ то
$f\in\Delta^{3} $ тоді і тільки тоді, коли $f'''(x)\ge 0,\ x\in (a,b).$

У статті йдеться про наближення $f\in\Delta^{3} $ кубічним сплайном, який теж з $\Delta^{3}$. А саме, ми доводимо
Теорему 1.

\medskip

\noindent {\bf Теорема 1.} {\it Нехай $\{a_j\}_{j=0}^n\ -$ набір рівновіддалених точок відрізку $[a,b],$ а саме
$\ a=a_n<a_{n-1}<...<a_0=b,\ h:=(b-a)/n.$ Якщо функція $f\in\Delta^{3},$ то існують
набір $\{b_j\}_{j=0}^n$ точок $[a,b]$ такий, що
$$
|a_j-b_j|\le 3h/2,\quad |b_j-b_{j-1}|\ge h/2,
$$
і
кубічний сплайн $s\in C^1$ з
 вузлами в точках $b_j$ такий, що
$$
s\in\Delta^{3} ,\eqno (1)
$$
$$
\left\Vert f-s \right\Vert_{[a_j,a_{j-1}]} \le c\, \omega_4
\left(f,h,[a_{j+4},a_{j-5}]\cap [a,b]\right), \quad
j=1,...,n, \eqno (2)
$$
а отже,
$$
\Vert f-s \Vert \le c\, \omega_4
\left(f,h,[a,b]\right), \eqno (3)
$$
де $\omega_4
\left(f,t,[\cdot,\cdot]\right)$ -- $4$-й модуль гладкості $f,\ a_\nu=a,\ \nu>n,$ та $a_\nu=b,\ \nu<0.$
}

\medskip

Тут і надалі $c$ позначають додатні абсолютні сталі, що можуть бути різними, навіть, якщо вони стоять в одному рядку.

Теорема 1 є частинним випадком \cite{6}, де, зокрема, побудовано і сплайн для оцінки з модулем гладкості Дітціана-Тотіка 4-го порядку.
В статті, пропонується простіша, ніж у \cite{6}, конструкція сплайна. Він,
на відміну від \cite{6}, представлений сумою усіченіх степеневих функцій, а отже, може бути використаний для побудови 3-монотонного многочлена, що наближатиме функцію як у (3).

Зауважимо, що навіть питання про справджуваність поточкового аналогу  (3),
тобто оцінки
$$
\vert f(x)-s(x) \vert \le c\, \omega_4
\left(f,\frac{1}{n^2}+\frac{\sqrt{1-x^2}}{n},[-1,1]\right), \quad x \in
[-1,1],
\eqno (4)
$$
є на сьогодні відкритим для 3-монотонного наближення (здається, що відповідь тут буде негативною, а функція $x^2\text{\rm sign} (x)\in\Delta^{3}$ буде контрприкладом, хоча ми не приділяємо цьому уваги у статті).
Також зауважимо, що неможливо замінити  $\omega_4$ на $\omega_k,\ k>4,$ у (3) (див. Шведов \cite{17}). Більш того, відповідні оцінки для 3-монотонного наближення у $L_p$ нормі з $p < \infty ,$
не є вірними навіть з  $\omega_3$ замість $\omega_4$
(див. Коновалов, Левіатан \cite{8} і Бондаренко, Примак \cite[Зауваження 5]{2}).

Що стосується $q$-монотонного наближення з $q>3,$ то тут (3) невірне також, навіть з
$\omega_3$ замість $\omega_4$ (див. \cite[Теорема 7.4]{6}), хоча для $1$-монотонного і $2$-монотонного наближень відповідні оцінки
справджуються.
Таким чином, можна сказати, що оцінка (3) є "граничною" \ між позитивними  і негативними  випадками у
 формозберігаючому наближенні, що розглядається. Розгорнутий огляд тематики див. в роботі Копотун, Левіатан, Примак,
Шевчук \cite{10}.

\medskip

\noindent {\bf Зауваження 1.}
 Сплайн $s,$ з Теореми 1, інтерполює $f$ у $a$ і $b,$ але, взагалі кажучи, це не інтерполяційний
 (у своїх вузлах)
сплайн степеня 3. Грубо кажучи, лише "мала" \ частина його вузлів залежить від
$f,$ тоді як решта є точками розбиття $a_j.$
Крім того, точки $a_j$ можуть бути "майже рівновіддаленими"\hspace{-.03cm}, див. Зауваження 2 у кінці статті.

\vskip 0.2 cm

Історія 3-монотонного наближення сплайнами є наступною.
Позначимо
$$
I:=[-1,1],\quad \|\cdot\|:=\|\cdot\|_{I},\quad
\omega_4\left(f,t\right):=\omega_4\left(f,t,I\right),\quad
\rho_n(x):=\frac{1}{n^2}+\frac{\sqrt{1-x^2}}{n}.
$$

Нехай $[a,b]=I.$ Для $f\in\Delta^{3}\cap C^2,$ Коновалов і Левіатан \cite{7} були першими, хто побудував квадратичний сплайн
 $s_1\in\Delta^{3}$ з $n$
рівномірними вузлами такий, що
$$
\left\Vert f-s_1 \right\Vert \le c\,n^{-2} \omega_1
\left(f'',1/n\right).
$$

Примак \cite{13} для $f\in\Delta^{3}$ побудував квадратичний сплайн
$s_2\in\Delta^{3}$ з $n$ довільними фіксованими вузлами такий, що, зокрема, для рівномірних вузлів
$$
\left\Vert f-s_2 \right\Vert \le c\, \omega_3 \left(f,1/n\right).
$$

З результів Шевчука \cite{16}, Левіатана і Примака
 \cite{11,14} випливає існування двох сплайнів $s_3$ і
 $s_4$ з $\Delta^{3}\cap C^3,$  обидва степеня 4 з $n$ рівномірними вузлами, таких, що
$$
\left\Vert f-s_3 \right\Vert \le c\,n^{-3} \omega_2
\left(f''',1/n\right),\quad n> 4,\quad f\in\Delta^{3}\cap C^3,
$$
$$
\left\Vert f-s_4 \right\Vert \le c\,n^{-1} \omega_4
\left(f',1/n\right),\quad n> N(f),\quad f\in\Delta^{3} \ (\subset C^1),
$$
де $N(f)$ є сталою, що залежить від $f.$ Зазначимо, що остання нерівність є, взагалі кажучи, не вірною
для всіх $n>4.$

Нещодавно Бондаренко, Левіатан і Примак \cite{1} довели першу
(і, мабуть, остаточну за порядком наближення)
поточкову оцінку
$$
\vert f(x)-S(x) \vert \le c\, \omega_3
\left(f,\rho_n(x)\right), \quad x \in I,
$$
з $S,$ що є 3-монотонним квадратичним сплайном по $n$-му чебишевському розбитті, і з $S,$ що є 3-монотонним многочленом степеня $\le n$.
\bigskip


\centerline{\bf \S 2. Допоміжні факти.} \vglue
0.3cm

\noindent
{\bf 2.1.} Нехай $\{a_j\}_{j=0}^n$ є множиною з $n+1$ фіксованих
точок $a_j:\ a=a_n<a_{n-1}<...<a_1<a_0=b,\ \ n\in\Bbb N.$ Для кожного
$j=q,...,n,\  q\in\Bbb N,\ q\le n, $ нехай
$L_q(x;a_j;g):=L(x;a_j,...,a_{j-q};g)$ позначає многочлен Лагранжа степеня $\le q$, який інтерполює
 $g\in C$ у $a_j,...,a_{j-q}.$

Покладемо
$$
S_q:=S_q(x):=\left\{ \aligned
L(x;a_q,...,a_0;g),&\quad x\in [a_q,b], \\
L(x;a_j,...,a_{j-q};g),&\quad x\in [a_j,a_{j-1}),\ \
j=q+1,...,n.
\endaligned
\right .\eqno (5)
$$

Для кожного $j=q,...,n,$ позначимо
$$
\Psi_q(x,a_j):=\left\{ \begin{aligned}
0,&\quad\text{якщо}\quad x\le a_j, \\
\prod\limits^{j}_{k=j-q+1}(x-a_k),&\quad\text{якщо}\quad x> a_j,
\end{aligned}
\right .
\quad\quad
\Psi_q(x,a_{q-1}):\equiv 0.
$$

\noindent {\bf Твердження \cite{5}.} {\it Неперервний сплайн $S_q$ має на
$[a,b]$ наспупне представлення
$$
S_q(x)=L_q(x;a_n;g)+\sum\limits_{j=q}^{n-1}
[a_{j+1},a_j,...,a_{j-q};g](a_{j-q}-a_{j+1}) \Psi_q(x,a_j),
$$
або еквівалентно,
$$
S_q(x)=L_{q-1}(x;a_n;g)+\sum\limits_{j=q}^{n}
[a_j,a_{j-1},...,a_{j-q};g]\left(\Psi_q(x,a_j)-\Psi_q(x,a_{j-1})\right),
$$
де квадратні дужки позначають розділені різниці $g.$
}

Для спрощення і не зменьшуючи загальності, далі будемо розглядати $[a,b]=I$.
Також, замість фігуруючих вище довільних фіксованих точок $\{a_j\}_{j=0}^n$ і замість
рівновіддалених точок $\{a_j\}_{j=0}^n$ у Теоремі 1, надалі візьмемо рівновіддалені, або
 "майже рівновіддалені"\ точки
$\{x_j\}_{j=0}^n,\ n\ge 3,\  -1=x_n<x_{n-1}<...<x_0=1.$ Позначимо
$$
I_j:=I_{j,n}:=[x_{j,n},x_{j-1,n}],\quad
h_j:=h_{j,n}:=x_{j-1,n}-x_{j,n},\quad j=1,...,n.
$$

Для $a\in I,$ покладемо
$$
\chi(x,a):=\left\{ \begin{aligned}
0,\quad&\text{якщо}\quad x\le a,\\
1,\quad&\text{якщо}\quad x>a,
\end{aligned}
\right . \quad (x-a)_+^r:=(x-a)^r\,\chi(x,a),\quad r\in\Bbb N.
$$
Будемо використовувати Твердження \cite{5} тільки для $q=3$ і $\{x_j\}_{j=0}^n$.
Тобто
$$
\Psi_3(x,x_{j})=(x-x_j)(x-x_{j-1})(x-x_{j-2})\chi(x,x_j)
$$
і
$$
S_3(x)=L_3(x;x_n;g)+\sum\limits_{j=3}^{n-1}
[x_{j+1},x_j,x_{j-1},x_{j-2},x_{j-3};g](x_{j-3}-x_{j+1})
\Psi_3(x,x_j),\eqno (6)
$$
або еквівалентно,
$$
S_3(x)=L_2(x;x_n;g)+\sum\limits_{j=3}^{n}
[x_j,x_{j-1},x_{j-2},x_{j-3};g]\left(\Psi_3(x,x_j)-\Psi_3(x,x_{j-1})\right).
\eqno(7)
$$
Нагадаємо, якщо $g\in\Delta^{3} ,$ то $[a,b,c,d;g]\ge 0$ для будь-яких різних
$a,b,c$ і $d.$ Зауважимо, що $S_3\notin\Delta^{3}$ навіть якщо $g\in\Delta^{3}.$

Будемо використовувати без спеціальних посилань нерівність Уітні
$
\left\Vert g-l_3\right\Vert_{[a,b]} \le \omega_{4}
\left(g,(b-a)/4,[a,b]\right),\ a<b,
$
де $l_3$ -- многочлен Лагранжа, що інтерполює $g$ у $a,\
a+\frac{b-a}{3},\ b-\frac{b-a}{3}$ і $b.$ Зразу відзначимо, що нерівності
$$
\left\Vert g-S_3\right\Vert_{I_j}\le c\,\omega_4\left(g,h_j,[x_j,x_{j-3}]\right),\quad j=1,...,n,
\eqno (8)
$$
випливають з (5) ($x_{-1}:=x_{-2}:=1$).

\medskip

{\bf 2.2.} Доведемо допоміжну Лему 1, яка складає і самостійний інтерес, якщо розглядається функція  $f$ з невід'ємними розділеними різницями порядку $q,\ q\ge 3.$ Для спрощення сформулюємо
і доведемо Лему 1 для $q=3.$ Позначимо
$$
\delta_j:=\delta_j(f):=[x_{j+1},x_j,x_{j-1},x_{j-2},x_{j-3};f], \quad
j=3,...,n-1.
$$
$$
\Delta_j:=\Delta_j(f):=[x_j,x_{j-1},x_{j-2},x_{j-3};f], \quad
j=3,...,n.
$$
Для спрощення позначень у Лемі 1, обмежимо $j$ значеннями $\{5,4,3\}$
і нехай $x_5<x_4<...<x_0$ -- будь-які фіксовані точки з $\{x_j\}_{j=0}^n.$

\bigskip

\noindent {\bf Lemma 1. } {\it Якщо $f\in\Delta^{3} ,$ то
$$
(x_1-x_4)(x_2-x_3)\Delta_4 \le (x_2-x_5)(x_3-x_4)\Delta_5
+(x_0-x_3)(x_1-x_2)\Delta_3 +
$$
$$
+2\left|\sqrt{(x_2-x_5)(x_2-x_4)\Delta_5(x_0-x_3)(x_1-x_3)\Delta_3}\right|=:A+2B.
\eqno (9)
$$
Більш того, якщо $\Delta_5\le\Delta_4\ge\Delta_3,$ то
$$
(x_1-x_4)(x_2-x_3)\Delta_4 \ge
\max\bigl\{(x_0-x_3)(x_1-x_2)\Delta_3-(x_2-x_5)(x_2-x_4+x_2-x_3)\Delta_5,\
$$
$$
(x_2-x_5)(x_3-x_4)\Delta_5-(x_0-x_3)(x_2-x_3+x_1-x_3)\Delta_3\bigl\}=:\max\{C,D\}\ge
A-2B.\eqno (10)
$$
}

\medskip

\noindent{\it Доведення.}  Використовуючи одно з представлень розділених різниць \cite{12},
(див. також у \cite[с.14]{15}), ми, для фіксованого $y\in(x_3,x_2),$ запишемо
$$
\begin{aligned}
\Delta_4(f)&=(x_2-x_4)[x_4,x_3,y,x_2;f]\Delta_4\bigl((x-x_3)(x-y)\chi(x,x_3)\bigl)+
\\
&+
(x_1-x_3)[x_3,y,x_2,x_1;f]\Delta_4\bigl((x-y)(x-x_2)\chi(x,y)\bigl)=
\\
&=\frac{y-x_4}{x_1-x_4}[x_4,x_3,y,x_2;f]+
\frac{x_1-y}{x_1-x_4}[x_3,y,x_2,x_1;f],
\\
\Delta_5(f)&=\frac{y-x_5}{x_2-x_5}[x_5,x_4,x_3,y;f]+
\frac{x_2-y}{x_2-x_5}[x_4,x_3,y,x_2;f],
\\
\Delta_3(f)&=\frac{y-x_3}{x_0-x_3}[x_3,y,x_2,x_1;f]+
\frac{x_0-y}{x_0-x_3}[y,x_2,x_1,x_0;f].
\end{aligned}
$$
Отже,
$$\begin{aligned}
\Delta_4(f)&=\frac{(y-x_4)(x_2-x_5)}{(x_2-y)(x_1-x_4)}\Delta_5(f)-
\frac{(y-x_4)(y-x_5)}{(x_2-y)(x_1-x_4)}[x_5,x_4,x_3,y;f]+
\\
&+\frac{(x_1-y)(x_0-x_3)}{(y-x_3)(x_1-x_4)}\Delta_3(f)-
\frac{(x_1-y)(x_0-y)}{(y-x_3)(x_1-x_4)}[y,x_2,x_1,x_0;f].
\end{aligned}
\eqno(11)
$$
Оскільки $y\in(x_3,x_2)$ і $f\in\Delta^{3}$ ($[a,b,c,d;f]\ge 0$), то маємо
$$
\Delta_4\le\min_{x_3<y<x_2}\left\{
\frac{(y-x_4)(x_2-x_5)}{(x_2-y)(x_1-x_4)}\Delta_5+
\frac{(x_1-y)(x_0-x_3)}{(y-x_3)(x_1-x_4)}\Delta_3\right\}.
$$
Знайшовши цей мінімум, бачимо, що $y_{\min}\in[x_3,x_2]$ для будь-яких
$\Delta_5,\Delta_3\ge 0$ і тому ми пишемо (9), беручи до уваги збіжність розділених різниць.

Доведемо першу нерівність у (10). Нехай $\max\{C,D\}=C.$ Припустимо обернене, що
$$
(x_1-x_4)(x_2-x_3)\Delta_4<C.\eqno(12)
$$
З (11) маємо
$$
\frac{(x_2-y)(x_1-x_4)}{(y-x_4)(x_2-x_5)}\left(\Delta_4-
\frac{(x_1-y)(x_0-x_3)}{(y-x_3)(x_1-x_4)}\Delta_3\right)\le\Delta_5,\quad
y\in(x_3,x_2).
$$
Разом з (12) це породжує
$$
E_1\Delta_4:=\left((x_1-x_4)(x_2-x_3)+(x_2-x_4+x_2-x_3)\frac{(x_2-y)(x_1-x_4)}{y-x_4}\right)\Delta_4<
$$
$$
\left((x_0-x_3)(x_1-x_2)+(x_2-x_4+x_2-x_3)\frac{(x_2-y)(x_1-y)(x_0-x_3)}
{(y-x_4)(y-x_3)(x_1-x_4)}\right)\Delta_3=E_2\Delta_3.
$$
Оскільки $\{x_j\}_{j=0}^n$
є рівновіддалені,  або "майже рівновіддалені" \,, тобто такі, що існує
$[a,b]\subset (x_3,x_2),$ для якого $E_1\ge E_2 > 0$ з $y\in[a,b],$ то остання
нерівність протирічить нерівності $\Delta_4\ge\Delta_3$ (зазначимо, що $a$ і $b,$ як нулі деякої параболи, можуть бути росташовані дуже близько один до одного і у найгіршому випадку порушення рівномірності точок $x_j,$ $a=b$). Випадок $C<D$ доводиться аналогічно.

Друга нерівність у (10) очевидна. Дійсно, якщо $C\ge D,$ то
$(x_0-x_3)(x_1-x_3)\Delta_3\ge (x_2-x_5)(x_2-x_4)\Delta_5$ і тому
$A-2B\le A-2|\sqrt{((x_2-x_5)(x_2-x_4)\Delta_5)^2}|=C.$
Лему 1 доведено.

\medskip

{\bf 2.3.} Зафіксуємо $n>3, \ j=3,...,n-1$ і будь-які
$a,c,b\in[x_{j+3},x_{j-5}]\cap I,\ a<c<b.$ Позначимо
$$
\widehat h_1:=c-a,\quad \widehat h_2:=b-c,\quad \widetilde
h_1:=b-x_j+b-x_{j-1}+b-x_{j-2},
$$
$$
\widetilde
h_2:=(b-x_j)(b-x_{j-1})+(b-x_j)(b-x_{j-2})+(b-x_{j-1})(b-x_{j-2}),$$$$
\widetilde h_3:=(b-x_j)(b-x_{j-1})(b-x_{j-2}).
$$
Означимо функцію  $\varphi_j\in C^1,$  співпадаючу з
$\Psi_3(x,x_j)$ майже скрізь,
$$
\varphi_j:=\varphi_j(x):=\varphi_j(x,a,c,b):=\alpha_j(x-a)_+^3+\beta_j(x-c)_+^3+\gamma_j(x-c)_+^2+(1-\alpha_j-\beta_j)(x-b)_+^3,
$$
де
$$
\alpha_j=\frac{\widetilde h_1\widehat h_2^2-2\widetilde h_2\widehat h_2+3\widetilde h_3}{3\widehat
h_1^2(\widehat h_1+\widehat h_2)},
$$
$$
\beta_j=\frac{\widetilde h_1\widehat h_2(\widehat h_1+\widehat h_2)(2\widehat h_1-\widehat
h_2)-\widetilde h_2(\widehat h_1^2-2\widehat h_2^2+2\widehat h_1\widehat h_2)-3\widetilde h_3(\widehat
h_2-\widehat h_1)}{3\widehat h_1^2\widehat h_2^2)},
$$
$$
\gamma_j:=\widetilde h_1-3\alpha_j(\widehat h_1+\widehat h_2)-3\beta_j\widehat h_2.
$$
Зауважимо, що
$$
\varphi_{j}(x)=\int\limits_{-1}^x\int\limits_{-1}^t \varphi_j''(u)dudt
.\eqno (13)
$$

\noindent {\bf Коментар.}  Числа $\alpha_j$ і $\beta_j$ обрані з двох відповідних умов так, щоб мати
(13) (перша робить рівним нулю добуток всіх коефіцієнтів біля $(x-\cdot )^1$ у представленні
 $\Psi_3$ сумою $(x-\cdot )^r,\ r=0,1,2,3,$ а друга  робить те саме з усіма
 вільними коефіцієнтами, включаючи той, що утворюється першою умовою).

\vskip 0.2 cm

\noindent Таким чином,
$$
\varphi_j(x,a,c,b)=\Psi_3(x,x_j),\quad x\in I\setminus
[\min\{a,x_j\},b]=:I\setminus \widehat {I}_j,\eqno(14)
$$
і якщо $h_j\le\widehat h_1<10h_j$ і $h_j\le\widehat h_2<10h_j,$ то
$$
\left\|\Psi_3(\cdot ,x_j)-\varphi_j\right\|=\left\|\Psi_3(\cdot
,x_j)-\varphi_j\right\|_{\widehat I_j}\le c\,h_j^3.\eqno(15)
$$

\bigskip
\medskip

\vspace{0.3cm}

\newpage

\centerline{\bf \S 3. Доведення Теореми 1.} \vglue
0.3cm

\centerline{\it Конструкція кубічного 3-монотонного сплайну}

\vskip 0.3 cm

\noindent {\bf 3.1.} Скрізь надалі $f\in\Delta^{3} .$ Нехай $n>4.$ Для кожного
$j=4,...,n-1$ позначимо $\Lambda_j:=(x_{j-3}-x_j)\Delta_j$ і якщо
$$
\Delta_{j+1}\le\Delta_j >\Delta_{j-1},\eqno(16)
$$
то будемо писати $j\in W.$ Для кожного $j\in W$ позначимо точку
$$
d_j:=\frac{(x_j+x_{j-1})\Lambda_{j+1}+(x_{j-1}+x_{j-2})\Lambda_{j}+(x_{j-2}+x_{j-3})\Lambda_{j-1}}
{2(\Lambda_{j+1}+\Lambda_j+\Lambda_{j-1})},
$$
що є центром параболи
$$
P_j(x):=\Lambda_{j+1}(x-x_{j})(x-x_{j-1})+\Lambda_{j}(x-x_{j-1})(x-x_{j-2})+
\Lambda_{j-1}(x-x_{j-2})(x-x_{j-3})=
$$
$$
=\bigl(\Lambda_{j+1}+\Lambda_{j}+\Lambda_{j-1}\bigl)(x-d_{j})^2+H_j,
$$
де
$H_j:=\Lambda_{j+1}(d_j-x_{j})(d_j-x_{j-1})+\Lambda_{j}(d_j-x_{j-1})(d_j-x_{j-2})+
\Lambda_{j-1}(d_j-x_{j-2})(d_j-x_{j-3}).$
Беручи до уваги першу нерівність у (10), нехай "майже рівновіддалені"\  $x_j$
є такі, що
$$
d_j\in I_{j-1},\quad j\in W\eqno(17)
$$
(для рівновіддалених $x_j$ таке включення гарановано).
Введемо
$$
Z:=\left\{j-1,j-2:\ j\in W\right\}
$$
(тобто, всі пари індексів, що відповідають кінцям проміжку у
 (17)). Використовуючи (9) і (10), легко перевірити, що
$$
\overline H_j\ge H_j\ge 0,\quad j\in W,\eqno(18)
$$
де $\overline H_j$ є найбільшим значенням $H_j,$ коли
$\Lambda_{j}(x_{j-2}-x_{j-1})=\Lambda_{j+1}(x_{j-1}-x_{j})+\Lambda_{j-1}(x_{j-3}-x_{j-2}).$ Аналогічно,
якщо у (9) ми маємо рівність, то $H_j=0.$

Нехай
$$
D:=\left\{d_j:\ j\in W\right\}.
$$
Відзначимо, що точки з $D$ розташовані на розбитті
 $x_j$ при наймні через один інтервал. Іншими словами, для будь-якого $j\in
W,$ індекси $j\pm 1,$ що відповідають розбиттю $x_j,$ не належать $W$ (у гіршому випадку тільки $j\pm 2$ можуть бути у
$W$). Зокрема, беручи до уваги (17), зауважимо, що якщо будь-яке
$$
j\in\left\{3,4,...,n-1\right\}=:J
$$
є таким, що $\Delta_{j+1}\le \Delta_j,$ то завжди $(x_j,x_{j-1})\cap
D=\emptyset ,$ а якщо воно таке, що $\Delta_{j+1}> \Delta_j,$ то завжди $(x_{j-1},x_{j-2})\cap D=\emptyset
.$ Останні два зауваження будуть використовуватися без спеціальних посилань.

Позначимо
$$
V:=J\setminus\bigl(W\cup\left\{j-1:\ j\in W\right\}\bigl).
$$
Замітимо, що $V=\left\{j:\ j-1\in J\setminus Z\right\}.$

Введемо
$$
Y:=\left\{y_i\right\}_{i=0}^k:=\left\{x_j:\ j\in (J\cup\{1,2\})\setminus Z
\right\}\cup D\cup \{-1,1\},
$$
де точки $y_i$ перенумеровано у зворотньому порядку і
$n-[n/3]-1\le k\le n.$

Далі, для кожного $j\in J,$ введемо нову функцію
$\Psi_j=\Psi_j(x)\in C^1.$
Для кожного $j\in V,$ через $i(j)$ позначимо такий індекс $i,$ при якому
$y_i=x_{j-1},$ і покладемо
$$
\Psi_j(x):=\left\{ \begin{aligned}
\varphi_{j}(x,y_{i(j)},y_{i(j)-1},y_{i(j)-2}),&\quad\text{якщо}\quad
\Delta_{j+1}\le \Delta_j,
\quad\quad\quad\quad\quad\quad\quad\quad\quad\quad\ \ \, \hfill (19)
\\
\varphi_{j}(x,y_{i(j)+2},y_{i(j)+1},y_{i(j)}),&\quad\text{інакше.}
\quad\quad\quad\quad\quad\quad\quad\quad\quad\quad\quad\quad\quad\quad\quad\quad\hfill (20)
\end{aligned}
\right .
$$
Для кожного $j\in W,$  через $i^*(j)$ позначимо такий індекс $i,$ при якому
$y_i=d_{j},$ і покладемо
$$
\Psi_j(x):=\varphi_{j}(x,y_{i^*(j)+1},y_{i^*(j)},y_{i^*(j)-1}),\eqno(21)
$$
$$
\Psi_{j-1}(x):=\varphi_{j-1}(x,y_{i^*(j)+1},y_{i^*(j)},y_{i^*(j)-1}).\eqno(22)
$$
Означення $\Psi_j,\ j\in J,$ завершено. Покладемо
$$
\Psi_{2}(x):=\Psi_{3}(x,x_2)=0,\quad
\Psi_{n}(x):=\Psi_{3}(x,x_n)=(x-x_n)(x-x_{n-1})(x-x_{n-2}).
$$
Таким чином, неперервно диференційовний на $I$ кубічний сплайн
$$
s(x):=L_3(x;x_n;f)+\sum\limits_{j=3}^{n-1}
\delta_j(f)\,(x_{j-3}-x_{j+1})\, \Psi_j(x),\eqno (23)
$$
або еквівалентно,
$$
s(x):=L_2(x;x_n;f)+\sum\limits_{j=3}^{n}
\Delta_j(f)\,\left(\Psi_j(x)-\Psi_{j-1}(x)\right),\eqno (24)
$$
має свої вузли лише в $Y$.

\vskip 0.3 cm

{\bf 3.2.} Доведем (3). Оскільки
$$
h_{j^*}\le|y_i-y_{i-1}|<4\,h_{j^*},\quad j=1,...,k,\eqno(25)
$$
де $h_{j^*}$ є довжина будь-якого найблищого до
$y_i$ проміжку $I_j$ (з двох можливих), то оцінка (3) випливає з  (8), (6), (23),
(14), (15) і оцінки
$$
\left|\delta_j\right|\le c \frac{\omega_4(f,h_j)}{h_j^4},\quad
j=3,...,n-1,
$$
див., наприклад, \cite[с.54]{15}. А саме, якщо $x\in I_{j^*},$  то
$$
\left|f(x)-s(x)\right|\le\left|f(x)-S_3(x)\right|+ \left|S_3(x)-s(x)\right|\le
$$
$$
\le c\,\omega_4(f,1/n)
+\sum\limits_{j=3}^{n-1}|\delta_j|(x_{j-3}-x_{j+1})\left|\Psi_3(x,x_j)-
\Psi_j(x)\right|=
$$
$$
=c\,\omega_4(f,1/n)
+\sum\limits_{j=\max\{3,j^*-5\}}^{\min\{n-1,j^*+4\}}
|\delta_j|(x_{j-3}-x_{j+1})\left|\Psi_3(x,x_j)- \Psi_j(x)\right|\le
c\,\omega_4(f,1/n).
$$
Отже, оцінку (3) доведено с одночасним доведенням (2).

Покажимо (1), тобто перевіремо, що $s''(x)$ неспадає на $I.$
 Заувжимо, що в усіх $\Psi_j,$
$$
\text{\rm sign}\gamma_j=\left\{ \begin{aligned} 1,&\quad\text{якщо}\quad
\Delta_{j+1}\le \Delta_j, \\
-1,&\quad\text{інакше,}
\end{aligned}
\right .\quad\quad j\in J. \eqno(26)
$$
Це буде зручно побачити разом з наступним.
Нехай $j\in V^+\ (V^-),$ якщо справджується (19) ((20)), відповідно, тобто $V=V^+\cup
V^-.$ Беручи до уваги  (25), замітимо, що у $\Psi_j$ з $j\in
V^+,\ \ \alpha_j+\beta_j\ge 0$ завдяки рівновіддаленому, або "майже рівновіддаленому"\
розбиттю $x_j,$ тоді як $\alpha_j\ge 0$ завжди
(для всіх $x_j$). Обидва числа (тобто, $\alpha_j+\beta_j$ і
$\alpha_j$) $\le 1.$ Аналогічно, якщо $j\in V^-,$ то $\alpha_j\le 1$ завдяки
$x_j,$  тоді як $\alpha_j+\beta_j\le 1$ завжди. Обидва числа $\ge
0.$ Будемо використовувати ці чотири зауваження зі спеціальним посиланням
 $(*).$

 \medskip

\noindent {\bf Зауваження 2.} Саме ці симетричні умови, разом з
(17) і Лемою 1, накладають обмеження на розбиття $x_j.$
Для рівновіддалених точок вони гарантовано виконуються, на відміну, скажімо, від чєбишевського розбиття
біля кінців інтервалу (у центральній частині чєбишевського розбиття теж все добре).

\vskip 0.2 cm

\noindent Стосовно чисел $\alpha_j$ і $\alpha_j+\beta_j$ у (21) а також
$\alpha_{j-1}$ і $\alpha_{j-1}+\beta_{j-1}$ у (22), замітемо, що
вони суттєво залежать тільки від розташування
центрального вузла
$y_{i^*(j)}=d_j.$ А саме, якщо $d_j$ знаходиться біля правого кінця
$I_{j-1}$ (див. (17)), то $\alpha_j$ і $\alpha_j+\beta_j$ "хороші"\,, тобто
задовольняють $(*)$ з $j\in V^+,$ тоді як $1\le
\alpha_{j-1}+\beta_{j-1}<1.5$ "погані"\  (тобто
не задовольняють $(*)$) і $0\le \alpha_{j-1}<0.5,$
інакше (якщо біля лівого кінця) $\alpha_{j-1}$ і
$\alpha_{j-1}+\beta_{j-1}$ "хороші"\,, тоді як $-0.5<\alpha_j\le 0$ і $1\le
\alpha_j+\beta_j<1.5$ "погані". (Фактично, $1.5$ і $-0.5$ це грубі числа.)
Будемо посилатись на "погані"\  властивості через $(**)$.

Таким чином, з (24) і (26) зразу помітимо, що
$$
s''(y_i-)\le s''(y_{i}+),\quad y_i\in Y.\eqno(27)
$$

Далі зазначимо, що суму у (24) зручно розглядати у спадному порядку, тобто від $n$ до $3,$
 як і дивитися на неспадність $s''$ а ні на невід'ємність $s'''$
на кожному $(y_i,y_{i-1}).$ Переконаємося, що
$$
s''(x)\nearrow ,\quad x\in (y_i,y_{i-1}),\ i=1,...,k.\eqno (28)
$$
Нехай $a_j\ (b_j)$ позначає найменший (найбільший) вузол з трьох вузлів кожної
 $\Psi_j,$ відповідно. Замітемо, що
$$
\Psi_j''(x)-\Psi_{j-1}''(x)=\left\{ \begin{aligned} 0,&\quad\text{якщо}\quad
x\in (-\infty, a_j], \\
2(x_{j-3}-x_j),&\quad\text{якщо}\quad x\in [b_{j-1},+\infty ),
\end{aligned}
\right .\quad\quad j=4,...,n. \eqno(29)
$$
Виділемо з (24) три доданки
$$
\Delta_{j+1}\bigl(\Psi_{j+1}''(x)-\Psi_{j}''(x)\bigl)+\Delta_{j}\bigl(\Psi_{j}''(x)-\Psi_{j-1}''(x)\bigl)+
\Delta_{j-1}\bigl(\Psi_{j-1}''(x)-\Psi_{j-2}''(x)\bigl)=:\overline P_j''(x),
\eqno (30)
$$
і розглянемо
$$
x\in (y_{i(j)+1},y_{i(j)})\cup(y_{i(j)},y_{i(j)-1})=:\widetilde
I_{i(j)+1}\cup \widetilde I_{i(j)}=:\overline{I}_{i(j)}.\eqno(31)
$$
Маємо три приципових ситуації: 1) якщо $j+1\in V^-,$ то $j\in
V;$ 2) якщо $j+1\in V^+,$ то $j\in V^+$ (і ніколи до $V^-$); 3)  $j+1\in
V\cup\{\nu -1:\ \nu\in W\},\ \ j\in W,\ \ j-2\in V\cup W.$ Перші два випадки є подібні, тому
перевіремо тільки другий.
Зауважимо, що $\Psi_{j+1}$ і $\Psi_{j}$ мають тільки два спільних вузла
$y_{i(j)}$ і $y_{i(j)-1}.$ Беручи до уваги (30) і (19),
запишемо
$$
s''(x)=A(x-y_{i(j)})+B=A(x-x_{j-1})+B\nearrow ,\quad x\in \widetilde
I_{i(j)},
$$
де $A\ge
(\Delta_{j+1}-\Delta_{j+2})(\alpha_{j+1}+\beta_{j+1})+(\Delta_{j}-\Delta_{j+1})\alpha_{j}\ge
0$ завдяки тому, що $\Delta_{j+2} \le \Delta_{j+1}\le \Delta_{j}$ разом з $(*),$
а $B$ є невід'ємна стала оскільки є (29).

Для отримання (28), у складному випадку 3), скористаємося тим, що
"хороші"\ властивості $(*)$ у сумі з "поганими"\  $(**)$ будь-що разом
 дадуть (28). Більш точно, беручи до уваги (14) і
(15), розглянемо для (31) три головні співвідношення (13), (18) і
(25). Нагадаємо, у цьому випадку $y_{i(j)}=y_{i^*(j)}=d_j.$ Завдяки (13)
маємо
$$
\int\limits_{-1}^x\int\limits_{-1}^t \overline P_j''(u)dudt=P_j(x)\chi
(x,d_j),\quad x\in I\setminus \overline I_{i^*(j)}, \eqno (32)
$$
і більш того,
$$
\int\limits_{-1}^x\ \overline P_j''(t)dt=P_j'(x)\chi (x,d_j),\quad x\in
I\setminus \overline I_{i^*(j)}. \eqno (33)
$$
Оскільки $\overline P_j''(x)$ має тільки три вузли і, більш того,
центральний це $d_j$ (тобто, центр $P_j$), то рівності (33)
і (32) з $H_j\ge 0$ не можуть бути вірними обидві разом одночасно без
(28), див. також (27). У якості додаткового зауваження, скажимо, що
нерівність (25) дає достатні відстані між цими трьома вузлами, щоб сформувати
 досить гарно обмежене число $\overline H_j$ (див. (18)) у (32) без того, щоб
  зруйнувати (28) на $\overline
I_{i^*(j)}$. Твердження (28), а отже і (1), доведено. Теорему 1 доведено.

\bigskip

\end{document}